# On an area property of the sum cotA + cotB + cotΓ in a triangle

## 1. *Introduction*

Rummaging through an obscure trigonometry book published in Athens, Greece (and in the Greek language), and long out of print, I discovered the following, listed as an exercise (see reference [1]):

**Quote:**

"Let $A\overset{\triangle}{B}\Gamma$ be a triangle. If on the side AB we draw the line perpendicular to AB at the point B; on the side BΓ the line perpendicular to it at the point Γ; and on ΓA the perpendicular to it at A; a new triangle A´B´Γ´ is formed. Prove that $\dfrac{E´}{E}$ = (cotA + cotB + cotΓ)$^2$, where E´ and E stand for the areas of the triangles $A´\overset{\triangle}{B´}\Gamma´$ and $A\overset{\triangle}{B}\Gamma$ respectively".

*end of quote*

# Important Note Regarding Angle Notation

In the interest of clarity and avoidance of confusion on the part of the reader, please note:

1) By refering to the angles A,B,Γ, we always be refering to the internal triangles A,B,Γ, in the context of a clearly defined triangle $A\overset{\triangle}{B}\Gamma$.

   In that context,
   angle A = < BAΓ = < ΓAB; or by using the alternative notation ^,

   the alternative notation **^,**

   angle A = $\overset{\wedge}{BA\Gamma}$ = $\overset{\wedge}{\Gamma AB}$.

   Similarly, angle B = < ABΓ = < ΓBA and

   angle Γ = $\overset{\wedge}{A\Gamma B}$ = $\overset{\wedge}{B\Gamma A}$.

   Consequently, all six angles A,B,Γ,A',B',Γ', referred to in Figures 1 through 4, are the internal angles of the triangles ABΓ and A'B'Γ'.

2) In any other context, we would use any of the two standard three-letter angle notations (mentioned above), with the vertex letter situated in the middle. However, in this paper, such "other context", does not truly arise.

3) In a paper like this, where the internal angles A,B,Γ are repeated very frequently throughout, using a three-letter notation (instead of a single letter one), would conceivably result in a significant increase in the length of the paper.

A first observation reveals that the sum cotA + cotB + cotΓ can be arbitrarily large. Indeed, if we restrict our attention only to those triangles in which none of the angles A, B, Γ is obtuse; then each of the trigonometric numbers cotA, cotB, cotΓ is bounded below by zero, i.e. cotA, cotB, cotΓ ≥ 0. But say, angle Γ can become arbitrarily small; in the language of calculus, Γ → $0^+$; so that cotΓ → + ∞. A second observation shows that the triangles $A\overset{\triangle}{B}\Gamma$ and $A´\overset{\triangle}{B}´\Gamma´$ are always similar (see Figures 2, 3, 4).

Furthermore note that if we were to effect a similar construction of triangle $A´\overset{\triangle}{B}´\Gamma´$ from a given triangle $A\overset{\triangle}{B}\Gamma$ but using an angle φ, 0° < φ ≤ 90° (instead of just φ = 90°), the same fact would emerge; that is, the two triangles are similar (see **Figure 1**).

This paper has a two-fold aim: first establish that $\frac{E´}{E}$ = (cotA + cotB + cotΓ)$^2$; and secondly that the minimum value of this area ratio is 3, attained precisely when the triangle $A\overset{\triangle}{B}\Gamma$ is equilateral (and thus, by similarity, triangle $A´\overset{\triangle}{B}´\Gamma´$ as well).

Even though, as we shall see, among **all** triangle pairs ($A\overset{\triangle}{B}\Gamma$, $A´\overset{\triangle}{B}´\Gamma´$), the ones with the smallest ratio (= 3) are those pairs in which both triangles are equilateral; if we restrict our search only to those pairs in which both $A\overset{\triangle}{B}\Gamma$ and $A´\overset{\triangle}{B}´\Gamma´$ are right triangles, then, as we will show, the minimum value of the above area ratio is equal to 4; obtained precisely when both right triangles are isosceles.

## 2. *Illustrations*

Angle Γ' = Angle A
Angle B' = Angle Γ
Angle A' = Angle B
$0° < A, B, Γ, A', B', Γ' < 90°$
$0° < φ < 90°$

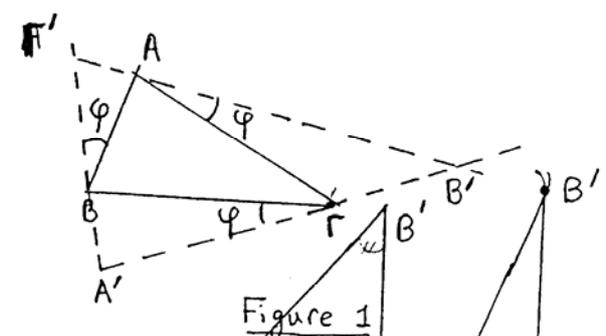
Figure 1

Angle Γ' = Angle A
Angle B' = Angle Γ
Angle A' = Angle B
$0° < A, B, Γ, A', B', Γ' < 90°$
$φ = 90°$

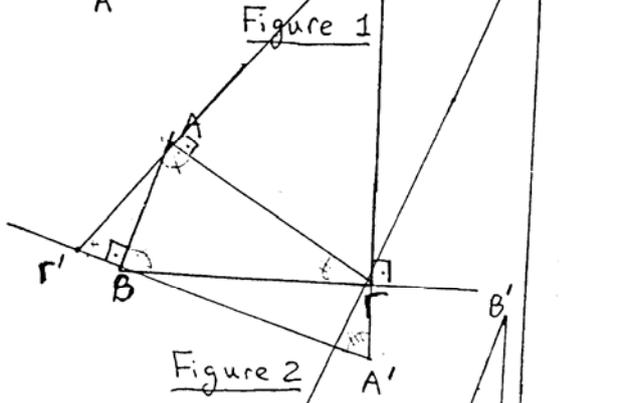
Figure 2

Point Γ' = Point B
Angle Γ' = Angle A = 90°
Angle B' = Angle Γ
Angle A' = Angle B
$0° < B, Γ, B', Γ' < 90°$
$φ = 90°$

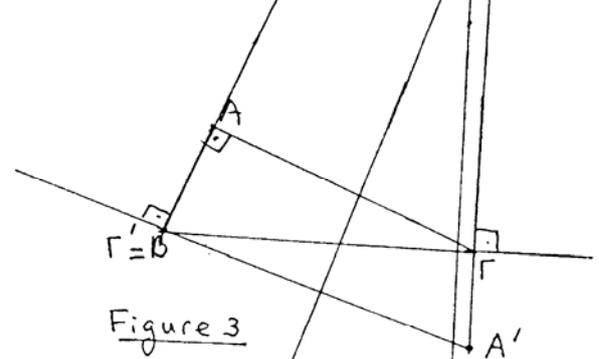
Figure 3

$0° < B, Γ < 90°$
$90° < A < 180°$
Angle Γ' = Angle A
Angle B' = Angle Γ
Angle A' = Angle B
$φ = 90°$

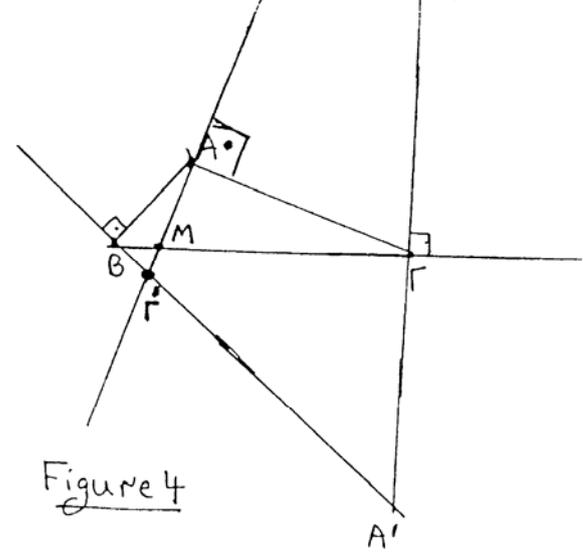
Figure 4

## 3. A few preliminaries

> The nine formulas listed below pertain to a triangle with sidelengths α = $|\overline{BΓ}|$, β = $|\overline{ΓA}|$, γ = $|\overline{AB}|$, angles A, B, Γ, semi-perimeter r (i.e. 2r = α + β + γ) and area E. Of these nine formulas, the first five are well known to a wide mathematical audience, while the last four (Formulas 6, 7, 8, 9) are less so; for each of those formulas, Formulas 6, 7, 8, and 9 we offer a short proof in the next section.
>
> **Formula 1** (Heron's Formula): E = $\sqrt{r(r-α)(r-β)(r-γ)}$
>
> **Formula 2:** cos2θ = 2cos²θ - 1 = 1 - 2sin²θ for any angle θ, typically measured in degrees or radians
>
> **Formula 3** (Law of Cosines): α² = β² + γ² - 2βγcosA, β² = α² + γ² - 2αγcosB, γ² = α² + β² - 2αβcosΓ
>
> **Formula 4** (summation formulas): For any angles θ and ω, cos (θ + ω) = cosθcosω - sinθsinω, sin (θ + ω) = sinθcosω + cosθsinω
>
> **Formula 5:** E = $\frac{1}{2}$ βγsinA = $\frac{1}{2}$ αβsinΓ = $\frac{1}{2}$ γαsinB
>
> **Formula 6:** 16E² = 2 (α²β² + β²γ² + γ²α²) - (α⁴ + β⁴ + γ⁴)
>
> **Formula 7:** For any angle θ which is not the form κπ or κπ + $\frac{π}{2}$, where κ is an integer, cot2θ = $\frac{\cot^2 θ - 1}{2\cot θ}$

**Formula 8:** $\cot\left(\dfrac{A}{2}\right) = \sqrt{\dfrac{r(r-\alpha)}{(r-\beta)(r-\gamma)}}$, $\cot\left(\dfrac{B}{2}\right) = \sqrt{\dfrac{r(r-\beta)}{(r-\alpha)(r-\gamma)}}$,

$\cot\left(\dfrac{\Gamma}{2}\right) = \sqrt{\dfrac{r(r-\gamma)}{(r-\beta)(r-\alpha)}}$

**Formula 9:** $E = \dfrac{\alpha^2 + \beta^2 + \gamma^2}{4(\cot A + \cot B + \cot \Gamma)}$

4. *Proofs of Formulas 6, 7, 8, and 9*

   a) **Proof of Formula 6:** From **Formula 1** we have,

   $E = \sqrt{r(r-\alpha)(r-\beta)(r-\gamma)} \Leftrightarrow 16E^2 = 16r(r-\alpha)(r-\beta)(r-\gamma);$

   $16E^2 = \left(\dfrac{\alpha+\beta+\gamma}{2}\right)\left(\dfrac{\beta+\gamma-\alpha}{2}\right)\left(\dfrac{\alpha+\gamma-\beta}{2}\right)\left(\dfrac{\alpha+\beta-\gamma}{2}\right) \cdot 16$

   $= (\alpha+\beta+\gamma)(\beta+\gamma-\alpha)[\alpha-(\beta-\gamma)](\alpha+\beta-\gamma)$

   $= [(\beta+\gamma)^2 - \alpha^2][\alpha^2 - (\beta-\gamma)^2]$

   $= \alpha^2(\beta+\gamma)^2 - \alpha^4 - [(\beta+\gamma)(\beta-\gamma)]^2 + \alpha^2(\beta-\gamma)^2$

   $= \alpha^2[(\beta+\gamma)^2 + (\beta-\gamma)^2] - \alpha^4 - [(\beta^2-\gamma^2)]^2$

   $= \alpha^2[2(\beta^2+\gamma^2)] - \alpha^4 - \beta^4 - \gamma^4 + 2\beta^2\gamma^2$

   $= 2(\alpha^2\beta^2 + \beta^2\gamma^2 + \gamma^2\alpha^2) - (\alpha^4 + \beta^4 + \gamma^4)$ ❑

   b) **Proof of Formula 7:** For any two angles θ and ω such that θ + ω ≠ κπ, κ any integer (so that sin (ω + θ) ≠ 0) we have, in accordance with **Formula 4**,

   $\cot(\theta + \omega) = \dfrac{\cos(\theta+\omega)}{\sin(\theta+\omega)} = \dfrac{\cos\theta\cos\omega - \sin\theta\sin\omega}{\sin\theta\cos\omega + \sin\omega\cos\theta}$

If, in addition to the above condition, we also have θ ≠ mπ and ω ≠ nπ (where m, n are any integers), then sinθsinω ≠ 0; and so,

$$\cos(\theta + \omega) = \frac{\frac{\cos\theta\cos\omega}{\sin\theta\sin\omega} - \frac{\sin\theta\sin\omega}{\sin\theta\sin\omega}}{\frac{\sin\theta\cos\omega}{\sin\theta\sin\omega} - \frac{\sin\omega\cos\theta}{\sin\omega\sin\theta}} = \frac{\cot\theta\cot\omega - 1}{\cot\omega + \cot\theta}$$

Thus, by setting θ = ω we obtain $\cot 2\theta = \frac{\cot^2\theta - 1}{2\cot\theta}$, under the conditions θ ≠ ℓπ and 2θ ≠ ℓπ; ℓ any integer. Considering the cases ℓ even and ℓ odd, the condition can be simplified into θ ≠ κπ and θ ≠ κπ + $\frac{\pi}{2}$, κ any integer. ❑

c) **Proof of Formula 8:** We need only establish the first of the three (sub)formulas, the other two are established similarly. Indeed, if we apply the first of the two (sub)formulas in **Formula 2** with θ = $\frac{A}{2}$; in combination with the first (sub)formula in **Formula 3** we obtain,

$$\cos^2\left(\frac{A}{2}\right) = \frac{1 + \frac{\beta^2 + \gamma^2 - \alpha^2}{2\beta\gamma}}{2} \Leftrightarrow \cos^2\left(\frac{A}{2}\right) = \frac{(\beta + \gamma)^2 - \alpha^2}{4\beta\gamma}$$

$$\Leftrightarrow \cos^2\left(\frac{A}{2}\right) = \frac{(\beta + \gamma - \alpha)(\beta + \gamma + \alpha)}{4\beta\gamma} = \frac{2(r - \alpha)(2r)}{4\beta\gamma} = \frac{(r - \alpha)r}{\beta\gamma} \quad (1)$$

Likewise, we apply $\cos 2\theta = 1 - 2\sin^2\theta$ with θ = $\frac{A}{2}$ in conjunction with the first (sub)formula in **Formula 3** to obtain,

$$\sin^2\left(\frac{A}{2}\right) = \frac{1 - \left(\frac{\beta^2 + \gamma^2 - \alpha^2}{2\beta\gamma}\right)}{2} = \frac{\alpha^2 - (\beta - \gamma)^2}{4\beta\gamma}$$

$$\Leftrightarrow \sin^2\left(\frac{A}{2}\right) = \frac{(\alpha + \gamma - \beta)(\alpha + \beta - \gamma)}{4\beta\gamma} = \frac{2(r-\beta)\,2(r-\gamma)}{4\beta\gamma}$$

$$= \frac{(r-\beta)(r-\gamma)}{\beta\gamma} \qquad (2)$$

Equations (1) and (2) imply

$$\cot^2\left(\frac{A}{2}\right) = \frac{\cos^2(A/2)}{\sin^2(A/2)} = \frac{r(r-\alpha)}{(r-\beta)(r-\gamma)} \qquad (3)$$

Since A is a triangle angle, $0° < \frac{A}{2} < 90°$; which means, that in particular $\cot\left(\frac{A}{2}\right) > 0$; thus, by virtue of (3) we conclude that

$$\cot\left(\frac{A}{2}\right) = \sqrt{\frac{r(r-\alpha)}{(r-\beta)(r-\gamma)}}. \qquad \square$$

d) **Proof of Formula 9:** We apply the first (sub)formula of **Formula 5**: $\beta\gamma = \frac{2E}{\sin A}$; we then substitute for $\beta\gamma$ in the first (sub)formula of **Formula 3** to obtain,

$$\alpha^2 = \beta^2 + \gamma^2 - \left(\frac{4E}{\sin A}\right)\cos A \Leftrightarrow \alpha^2 = \beta^2 + \gamma^2 - 4E\cot A \qquad (4)$$

Similarly, we also have, $\qquad \beta^2 = \alpha^2 + \gamma^2 - 4E\cot B \qquad (5)$

and, $\qquad \gamma^2 = \alpha^2 + \beta^2 - 4E\cot\Gamma \qquad (6)$

Adding equations (4), (5), and (6) memberwise and solving for the area E produces the desired result: $E = \frac{\alpha^2 + \beta^2 + \gamma^2}{4(\cot A + \cot B + \cot\Gamma)}. \qquad \square$

## 5. *Postulate 1 and its proof*

**Postulate 1**

> Let $A\overset{\Delta}{B}\Gamma$ be a triangle and $A'\overset{\Delta}{B'}\Gamma'$ be the triangle constructed from $A\overset{\Delta}{B}\Gamma$ by drawing the following three straight lines: The perpendicular to AB at the point B, the perpendicular to BΓ at Γ, and the perpendicular to ΓA at A. If E and E´ are the areas of $A\overset{\Delta}{B}\Gamma$ and $A'\overset{\Delta}{B'}\Gamma'$ respectively, $\dfrac{E'}{E} = (cotA + cotB + cot\Gamma)^2$.

**Proof:** Without loss of generality, we may assume that B, Γ are acute angles; 0° < B, Γ < 90°; while A maybe an acute, right, or obtuse angle. The key point in this proof is the observation or realization that,

$$E' = E + \frac{\gamma^2 cotA + \beta^2 cot\Gamma + \alpha^2 cotB}{2} \qquad (7)$$

To see why this is so let us first take a look at **Figure 2**. In that picture we see that $A\overset{\Delta}{B}\Gamma$ lies in the interior of $A'\overset{\Delta}{B'}\Gamma'$ and that,

$$E' = E + (\text{area } A\overset{\Delta}{B}\Gamma') + (\text{area } B\overset{\Delta}{\Gamma}A') + (\text{area } \Gamma\overset{\Delta}{A}B'),$$

from which (7) easily follows (look at the formulas listed in the space adjacent the illustration in **Figure 2**; and just apply the area formula for a right triangle).

Now, in the case of A = 90° (**Figure 3**), (area $A\overset{\Delta}{\Gamma'}B$) = 0 (since cotA = cot90° = 0); equation (7) still holds true with the term γ²cotA being zero.

Next, let us look at **Figure 4**. Here, triangle $A\overset{\triangle}{B}\Gamma$ does not lie in the interior of $A'\overset{\triangle}{B'}\Gamma'$; however, as we show below, equation (7) still holds true; the piece of $A\overset{\triangle}{B}\Gamma$ that lies outside the interior of triangle $A'\overset{\triangle}{B'}\Gamma'$ (that is, triangle $A\overset{\triangle}{B}M$), is compensated by the fact that cotA < 0, on account of 90° < A < 180°. Indeed,

$E'$ = (area $B'\overset{\triangle}{A}\Gamma$) + (area $A\overset{\triangle}{M}\Gamma$) + (area of the quadrilateral MΓA'Γ');

$E'$ = (area $B'\overset{\triangle}{A}\Gamma$) + E - (area $A\overset{\triangle}{M}B$) + (area $B\overset{\triangle}{\Gamma}A'$) - (area $B\overset{\triangle}{M}\Gamma'$);

and by virtue of (area $A\overset{\triangle}{B}\Gamma'$) = (area $A\overset{\triangle}{M}B$) + (area $B\overset{\triangle}{M}\Gamma'$),

we obtain $E'$ = E + (area $B'\overset{\triangle}{A}\Gamma$) + (area $B\overset{\triangle}{\Gamma}A'$) - (area $A\overset{\triangle}{B}\Gamma'$)     (8)

We have, (area $B'\overset{\triangle}{A}\Gamma$) = $\dfrac{\beta^2 \cot\Gamma}{2}$, (area $B\overset{\triangle}{\Gamma}A'$) = $\dfrac{\alpha^2 \cot B}{2}$, and

(area $A\overset{\triangle}{B}\Gamma'$) = $\dfrac{\gamma \cdot |\overline{B\Gamma'}|}{2}$; and since (from triangle $A\overset{\triangle}{B}\Gamma'$ in **Figure 4**) we have $|\overline{B\Gamma'}|$ = γcot (180° - A) = - γcotA, we conclude that (area $A\overset{\triangle}{B}\Gamma'$)

= - $\dfrac{\gamma^2 \cot A}{2}$; and by (8) we arrive at equation (7); thus (7) holds true in all cases. Next note that,

$$\dfrac{E'}{E} = (\cot A + \cot B + \cot\Gamma)^2 \Leftrightarrow \quad \text{(by \textbf{Formula 9})}$$

$$\Leftrightarrow \quad \dfrac{E'}{E} = \left(\dfrac{\alpha^2 + \beta^2 + \gamma^2}{4E}\right)^2 \Leftrightarrow 16EE' = (\alpha^2 + \beta^2 + \gamma^2)^2 \Leftrightarrow$$

(by (7)) $\Leftrightarrow \quad 16E\left(E + \dfrac{\gamma^2 \cot A + \beta^2 \cot\Gamma + \alpha^2 \cot B}{2}\right) = (\alpha^2 + \beta^2 + \gamma^2)^2 \Leftrightarrow$

$$\Leftrightarrow \quad 16E^2 + 8\,(\gamma^2\cot A + \beta^2\cot\Gamma + \alpha^2\cot B)E - (\alpha^2 + \beta^2 + \gamma^2)^2 = 0 \quad (9)$$

Clearly, if we prove (9), the proof of **Postulate 1** will be complete.

First note that according to **Formula 6** we have,

$$16E^2 = 2\,(\alpha^2\beta^2 + \beta^2\gamma^2 + \gamma^2\alpha^2) - (\alpha^4 + \beta^4 + \gamma^4) \quad (10)$$

Secondly, we compute the term $8E\gamma^2\cot A$ in (9), in terms of $\alpha$, $\beta$, $\gamma$; and by cyclicity of the letters $\alpha$, $\beta$, $\gamma$, we correspondingly find similar formulas for the terms $8E\beta^2\cot\Gamma$ and $8E\alpha^2\cot B$. We start by applying the first (sub)formula of **Formula 8** together with **Formulas 7** and **1** to obtain,

$$8E\gamma^2\cot A = 8\,\sqrt{r(r-\alpha)(r-\beta)(r-\gamma)} \cdot \frac{1}{2} \cdot \left[\frac{r(r-\alpha)}{(r-\beta)(r-\gamma)} - 1\right] \cdot \sqrt{\frac{(r-\beta)(r-\gamma)}{r(r-\alpha)}}$$

$$\Leftrightarrow 8E\gamma^2\cot A = 4\gamma^2\,[r(r-\alpha) - (r-\beta)(r-\gamma)] \Leftrightarrow$$

$$\Leftrightarrow 8E\gamma^2\cot A = 4\gamma^2\left[\left(\frac{\alpha+\beta+\gamma}{2}\right)\left(\frac{\beta+\gamma-\alpha}{2}\right) - \left(\frac{\alpha+\gamma-\beta}{2}\right)\left(\frac{\alpha+\beta-\gamma}{2}\right)\right] \Leftrightarrow$$

$$\Leftrightarrow 8E\gamma^2\cot A = \gamma^2\,[(\beta+\gamma)^2 - \alpha^2 - \{\alpha^2 - (\beta-\gamma)^2\}] \Leftrightarrow$$

$$\Leftrightarrow 8E\gamma^2\cot A = 2\gamma^2\,(\beta^2 + \gamma^2 - \alpha^2) \quad (11)$$

Similarly, by cyclicity, $8E\beta^2 \cot\Gamma = 2\beta^2 (\alpha^2 + \beta^2 - \gamma^2)$ and  (12)

$8E\alpha^2 \cot B = 2\alpha^2 (\gamma^2 + \alpha^2 - \beta^2)$  (13)

Also, by expansion we have,

$- (\alpha^2 + \beta^2 + \gamma^2)^2 = -2(\alpha^2\beta^2 + \beta^2\gamma^2 + \gamma^2\alpha^2) - (\alpha^4 + \beta^4 + \gamma^4)$  (14)

If we add (10), (11), (12), (13), and (14) memberwise, the resulting equation is (9). ❏

## 6. *Two Lemmas from calculus and their conclusion*

> **Lemma 1:** Over the open interval $(0, +\infty)$, the function $f(x) = \dfrac{2x^2 - x\sqrt{x^2+1} + 2}{\sqrt{x^2+1}}$ has an absolute minimum at $x = \dfrac{1}{\sqrt{3}}$; the absolute minimum value being $f\left(\dfrac{1}{\sqrt{3}}\right) = \sqrt{3}$.
>
> **Lemma 2:** Let k be a positive constant. Over the open interval $(0, \dfrac{\pi}{2})$, the function $g_k(x) = \dfrac{\cot^2 x + k\cot x + k^2 + 1}{\cot x + k}$ has an absolute minimum at $x = \theta_k$; where $\theta_k$ is the unique real number in $(0, \dfrac{\pi}{2})$ such that $\cot\theta_k = -k + \sqrt{k^2+1}$; the absolute minimum value being $g_k(\theta_k) = \dfrac{2k^2 - k\sqrt{k^2+1} + 2}{\sqrt{k^2+1}}$ (it can be easily verified that $g_k(\theta_k) > 0$.)

The proof to the first Lemma utilizes the standard calculus techniques and is left to the reader. The proof of the second Lemma can also be done by the use of the standard techniques.

To the interested reader we point out that the derivative of the function $g_k$ is given by

$$g'_k(x) = \frac{-\csc^2 x \cdot [\cot^2 x + 2k\cot x - 1]}{(\cot x + k)^2} = \frac{-\csc^2 x \cdot [\cot x + k - \sqrt{k^2+1}]}{(\cot x + k)^2}$$

$[\cot x + k + \sqrt{k^2 + 1}]$

Note that in the last fraction, only the factor $\cot x + k - \sqrt{k^2+1}$ changes sign over the interval $(0, \frac{\pi}{2})$; that happens at the only critical number $\theta_k$ this function has: the unique number $\theta_k$ such that $\cot \theta_k = -k + \sqrt{k^2+1}$ (recall that $\cot x$ is a decreasing function on $(0, \frac{\pi}{2})$ and whose range is $(0, +\infty)$).

Combining the two Lemmas we arrive at,

**The conclusion from the two Lemmas**

> *Among all functions $g_k$ in Lemma 2, the one with the smallest absolute minimum is the function $g_k$ with $k = \frac{1}{\sqrt{3}}$. The smallest absolute minimum occurs at*
>
> $$x = \theta_k = \frac{\pi}{3} \text{ (from } \cot \theta_{1/\sqrt{3}} = -\frac{1}{\sqrt{3}} + \sqrt{\left(\frac{1}{\sqrt{3}}\right)^2 + 1} = \frac{1}{\sqrt{3}}\text{)}.$$
>
> *The smallest minimum is given by* $g_{1/\sqrt{3}}\left(\frac{1}{\sqrt{3}}\right) = \sqrt{3}$.

## 7. *Postulates 2 and 3*

Now let us look at the sum cotA + cotB + cotΓ with 0° < B, Γ < 90° and 0° < A < 180°. Since A = 180° - (B + Γ) we have,

cotA + cotB + cotΓ = cot [180° - (B + Γ)] + cotB + cotΓ =

= - cot (B + Γ) + cotB + cotΓ.

Recall from the proof of **Formula 7** that $\cot(B + \Gamma) = \dfrac{\cot B \cot \Gamma - 1}{\cot B + \cot \Gamma}$.

Therefore cotA + cotB + cotΓ = $-\left(\dfrac{\cot B \cot \Gamma - 1}{\cot B + \cot \Gamma}\right)$ + cotB + cotΓ;

cotA + cotB + cotΓ = $\dfrac{\cot^2 B + \cot \Gamma \cot B + \cot^2 \Gamma + 1}{\cot B + \cot \Gamma}$.

If we hold the angle Γ fixed, cotΓ becomes a positive constant; cotΓ = k. And if we set B = x, and let x vary in the open interval (0, $\dfrac{\pi}{2}$), the above cotangent sum becomes simply one of the functions $g_k$ of **Lemma 2**! Hence, according to the **Conclusion of the two Lemmas**, the smallest value of the above sum is equal to $\sqrt{3}$ attained when k = cotΓ = $\dfrac{1}{\sqrt{3}}$ and B = x = $\dfrac{\pi}{3}$; in degrees, B = Γ = 60°; A = 60° as well.

Combining this with **Postulate 1** leads to

## Postulate 2

> *The ratio $\frac{E'}{E}$ has the minimum value 3. The minimum value is attained precisely among the triangle pairs $(A\overset{\triangle}{B}\Gamma, A'\overset{\triangle}{B'}\Gamma')$, with both triangles being equilateral.*

We finish by taking a look at the cases in which $A\overset{\triangle}{B}\Gamma$ (and thus $A'\overset{\triangle}{B'}\Gamma'$ as well) is a right triangle: A = 90°; Γ = 90° - B; so that

cotA + cotB + cotΓ = 0 + cotB + cot (90° - B) = cotB + tanB =

$$\frac{\cos B}{\sin B} + \frac{\sin B}{\cos B} = \frac{\cos^2 B + \sin^2 B}{\sin B \cos B} = \frac{1}{\sin B \cos B} = \frac{2}{2\sin B \cos B} = \frac{2}{\sin 2B} \geq 2,$$

since 0 < sin2B ≤ 1, in view of 0° < 2B < 180°. Thus the smallest value of (cotA + cotB + cotΓ)² is 4; this happens when sin2B = 1; B = Γ = 45°.

## Postulate 3

> *Among all pairs $(A\overset{\triangle}{B}\Gamma, A'\overset{\triangle}{B'}\Gamma')$ in which both triangles are right ones the minimum value of $\frac{E'}{E}$ is 4 attained when both triangles are isosceles.*